\documentclass[11pt]{article}

\usepackage{diagram,amssymb,latexsym}

\makeatletter

\def\section{\@startsection{section}{1}%
  {0pt}{1.37\baselineskip}{.5\baselineskip}%
  {\upshape\large\bfseries}}

\makeatother

\newtheorem{thm}[subsection]{Theorem}
\newtheorem{prop}[subsection]{Proposition}
\newtheorem{lemma}[subsection]{Lemma}
\newenvironment{prf}{{\noindent \hspace{1.3em} \bf Proof.}}
                     {\hspace*{\fill} $\blacksquare$ \vspace{4pt}}
\newtheorem{cor}[subsubsection]{Corollary}
\newenvironment{definition}{\vspace{.3\baselineskip} 
                            \noindent 
                            \hspace{-.75em} 
                            \addtocounter{subsection}{1}
                            {\bf Definition
                             \arabic{section}.\arabic{subsection}}}{}
\newenvironment{example}{\vspace{.5\baselineskip} 
                         \noindent 
                         \hspace{-.75em} 
                         \addtocounter{subsection}{1}
                          {\bf Example 
                          \arabic{section}.\arabic{subsection}}}{}

\newcommand{\ie}{{\it i.e.\ }}
\newcommand{\mut}{\mu_T}
\newcommand{\lambdal}{\lambda_L}
\newcommand{\deltal}{\delta_L}
\newcommand{\lambdar}{\lambda_R}
\newcommand{\deltar}{\delta_R}
\newcommand{\comod}{\bigtriangledown}

\newcommand{\Ae}{A^e}
\newcommand{\myrightarrow}[1]{\stackrel{#1}{\longrightarrow}}

\newcommand{\Hom}[1]
            {\mathrm{Hom}_{\Ae}\left(#1 \, , \, N \otimes M \right)}

\begin{document}

\begin{center}
   {\large \bf 
   Modules, comodules and cotensor products \\
   over Frobenius algebras
   }
   \\
   \medskip
   {\large Lowell Abrams} \\
   \vspace{1pt}
   {\small 
   Department of Mathematics \\
   Rutgers University \\
   New Brunswick, NJ 08903  \\
   \vspace{1pt}
   labrams@math.rutgers.edu
   }
\end{center}

\vspace{.5in}
	
\begin{abstract}
We characterize noncommutative Frobenius algebras $A$ in terms of the
existence of a coproduct which is a map of left $\Ae$-modules. We show
that the category of right (left) comodules over $A$, relative to this
coproduct, is isomorphic to the category of right (left) modules. This
isomorphism enables a reformulation of the cotensor product of
Eilenberg and Moore as a functor of modules rather than comodules.

We prove that the cotensor product $M \Box N$ of a right $A$-module
$M$ and a left $A$-module $N$ is isomorphic to the vector space of
homomorphisms from a particular left $\Ae$-module $D$ to $N \otimes
M$, viewed as a left $\Ae$-module. Some properties of $D$ are
described. Finally, we show that when $A$ is a symmetric algebra, the
cotensor product $M \Box N$ and its derived functors are given by the
Hochschild cohomology over $A$ of $N \otimes M$.

\vspace{1pt}
\noindent
Keywords: Frobenius algebra, comodule, cotensor product, Hochschild
cohomology
\end{abstract}

\section{Introduction}

Eilenberg and Moore originally introduced the cotensor product $M \Box
N$ and its derived functors Cotor$(M,N)$ on comodules $M, N$ as tools
for the calculation of the homology of the fiber space in a fibration
\cite{EilMoo66}. This paper investigates these functors in the context
where the coalgebra is a Frobenius algebra (defined in section
\ref{sec:FA}). 

The Frobenius case is not far removed from that of Eilenberg and
Moore, whose coalgebra is the set of normalized singular chains in
some space $X$; in the presence of sufficient flatness, all the
relevant constructions yield exactly the same data upon passing to
homology [{\it ibid.}]. When the space $X$ under consideration is
compact and oriented, its homology is in fact a Frobenius algebra.

Nevertheless, our approach diverges from that of Eilenberg and Moore
in an important way. The results presented here rest on a new
characterization of Frobenius algebras as algebras possessing a
coassociative counital comultiplication $\delta \colon A \rightarrow A
\otimes A$ which is a map of regular bimodules. (This is formulated
slightly differently as theorem \ref{th:comult} below.)  This
comultiplication is decidedly different from the one used by Eilenberg
and Moore. The relationship between the two coproducts will be
discussed elsewhere.

The Frobenius algebra coproduct, and in particular the element
$\delta(1_A)$, has already begun to find its place in a variety of
contexts. In two dimensional topological quantum field theory, it
gives rise to the handle operator \cite{Abr96}. In quantum cohomology
it provides a generalization of the classical Euler class
\cite{Abr97}. It also plays an important role in the study of quantum
Yang-Baxter equations and, under certain conditions, serves as a
separability idempotent \cite{BeiFonSto97}. Here, we will consider
left submodules of $A \otimes A$ generated by $\delta(1_A)$ and $T
\circ \delta(1_A)$. These will be discussed more later in this
section.

The bimodule property of the Frobenius algebra coproduct implies
another important property of Frobenius algebras, appearing as theorem
\ref{th:modcomod}: The category of right modules over a Frobenius
algebra $A$ is isomorphic to the category of right comodules over
$A$. This result makes it possible to view Eilenberg and Moore's
functors on comodules as functors on modules. Now, using the Snake
Lemma, one can show that the cotensor product is left exact in both
variables. (This also follows from theorem \ref{th:cotensorhom}, of
course.) This suggests that the right module $M \Box N$ should be
expressible as a module of homomorphisms from some left module $D$ to
$N \otimes M$. In fact, this is the case, as stated in theorem
\ref{th:cotensorhom}. The concern is to develop a satisfactory
understanding of the module $D$.

Specifically, $D$ denotes the left $\Ae$-submodule of $A \otimes A$
generated by $T \circ \delta(1_A)$, where $T \colon A \otimes A
\rightarrow A \otimes A$ denotes the canonical involution. This is not
the same as the left $\Ae$-submodule $\delta(A)$ of $A \otimes A$
generated by $\delta(1_A)$. The latter module is a very natural object
to consider, since $\delta$ itself is a left $\Ae$-module map, but the
importance of $D$ in this context is somewhat surprising. Under
certain conditions, delineated in \ref{pr:symmdelta} and
\ref{cor:Ddelta}, $D$ and $\delta(A)$ are in fact the same up to a
canonical involution. But in other cases, such as the one presented
below as example 4.\ref{ex:1}, this is not so.

There are two important corollaries to the main results discussed
above. One (\ref{cor:cotorext} below) is that the right derived
functors of the cotensor product $M \Box N$, \ie Cotor$^*(M,N)$, are
in fact the modules Ext$^*(D,N \otimes M)$. The other
(\ref{cor:cotorHoch} below) is that when $A$ is a symmetric algebra,
the cotensor product $M \Box N$ and its derived functors are given by
the Hochschild cohomology over $A$ of $N \otimes M$.

The author extends heartfelt thanks to Chuck Weibel who, in addition
to being free with helpful advice, is a living index to
\cite{Wei}. Peter May also offered several useful comments and
suggestions. 

\vspace{.25\baselineskip}

\noindent
{\bf Notation and Conventions}

All algebras $A$ considered here are assumed to be finite dimensional
as a vector space over their coefficient field $K$, and to possess a
multiplicative identity element $1_A$. We let $\mu \colon A \otimes A
\rightarrow A$ denote the multiplication map. The symbols $A^n$ will
always denote $A^{\otimes n}$, \ie the tensor product of $n$ copies of
$A$, and never the Cartesian product. For any object $X$, we will use
``X'' or ``$\cdot$'' to denote the identity map $X \rightarrow X$, and
the symbols $\cdot \otimes \cdot$ will be abbreviated ``$\cdot
\cdot$''.

\section{Noncommutative Frobenius Algebras} 
   \label{sec:FA}

An algebra $A$ is defined to be a {\bf Frobenius algebra} if it
possesses a left $A$-module isomorphism $\lambdal \colon A \rightarrow
A^*$ with its vector space dual.  Here, $A$ is viewed as the left
regular module over itself, and $A^*$ is made a left $A$-module by the
action $(a \cdot \zeta)(b) \, := \, \zeta(ba)$ for any $a,b \in A$ and
$\zeta \in A^*$. It is easy to show that the existence of the
isomorphism of left modules implies the existence of an isomorphism
$\lambdar$ of right modules, where the right module structures are
defined analogously.

There are many equivalent definitions of Frobenius algebras; see
\cite{CurRei} for more information. For our purposes, the new
characterization of Frobenius algebras presented below is very useful.

\begin{thm}
   \label{th:comult} 
An algebra $A$ is a Frobenius algebra if and only if it has a
coassociative counital comultiplication $\delta \colon A \rightarrow A
\otimes A$ which is a map of left $\Ae$-modules.
\end{thm}

Here, $\Ae$ denotes the ring $A \otimes A^{\mathrm{op}}$, and $A$ has
the left $\Ae$-action defined by $(b \otimes b') \cdot a:= bab'$.

In many respects, the proof of this result follows the proof of an
analogous result for the commutative case, found in \cite{Abr96}. For
the sake of space, we merely indicate how this proof differs from the
one given there.
                                         
\begin{prf}
Assume $A$ denotes a Frobenius algebra with left-module isomorphism
$\lambdal \colon A \rightarrow A^*$. Let $\mut \colon A \otimes A
\rightarrow A$ denote the composition $\mu \circ T$. Define the
comultiplication map $\deltal \colon A \rightarrow A \otimes A$ to be
the composition $(\lambdal^{-1} \otimes \lambdal^{-1}) \circ \mut^*
\circ \lambdal$. With the appropriate adjustments, the discussion in
\cite{Abr96} shows that the following diagram commutes:
\[
\begin{diagram}
   \node{A \otimes A} \arrow{e,t}{\mu} 
                      \arrow{s,r}{\cdot \otimes \deltal}
       \node{A} \arrow{s,r}{\deltal} \\
   \node{A \otimes A \otimes A} \arrow{e,t}{\mu \otimes \cdot}
      \node{A \otimes A}
\end{diagram}
\]
In words, $\deltal$ is a map of left $A$-modules. 

Using the right-module isomorphism $\lambdar \colon A \rightarrow
A^*$, it is an analogous exercise to define $\deltar$ and show that
this comultiplication map is a map of right modules. Let $\epsilon
\colon A \rightarrow K$ denote $\lambdar(1_A)$. Note that
$\lambdar(1_A) = \lambdal(1_A)$, and thus that $\epsilon$ serves as a
counit for both $\deltar$ and $\deltal$.

Now consider the following diagram:
\[
\begin{diagram}
   \node[2]{A} \arrow{se,t}{\deltar} \\
    \node{A^2} \arrow{se,b}{\deltar \otimes \deltal}
                        \arrow{e,t}{\deltar \otimes \cdot}
                        \arrow{ne,t}{\mu}
       \node{A^3}
                 \arrow{s,r}{\cdot \otimes A \otimes \deltal}
                  \arrow{e,b}{\cdot \otimes \mu}
         \node{A^2}
                   \arrow{s,r}{\deltal}
                   \arrow{se,t}{\cdot \cdot} \\
      \node[2]{A^4} 
                        \arrow{e,b}{\cdot \otimes \mu \otimes A}
      \node{A^3} 
                    \arrow{e,b}{\cdot \otimes \epsilon \otimes \cdot}
      \node{A^2}    
\end{diagram}
\]
This diagram commutes because of the properties of $\deltar$,
$\deltal$ and $\epsilon$ mentioned just above. It follows that
$\deltar \circ \mu$ is the same as the composition of maps from the
far left down and along the bottom row to the lower right-hand
corner. A corresponding diagram shows that $\deltal \circ \mu$ is also
the same as that composition, \ie $\deltar \circ \mu = \deltal \circ
\mu$. Since $A$ has an identity element, we see that $\deltar =
\deltal$.  Define $\delta := \deltar = \deltal$. We have just shown
that this map $\delta \colon A \rightarrow A \otimes A$ is a map of
bimodules, \ie is an $\Ae$-module map, and has a counit.

The remainder of the proof follows as in \cite{Abr96}.
\end{prf}

Throughout the sequel, $\delta$ and $\epsilon$ will denote the
comultiplication and counit respectively. Let $\delta(A)$ denote the
image of $\delta$.

\begin{cor}
   \label{cor:delta}
The map $\delta$ is an injection of left $\Ae$-modules.
\end{cor}
\begin{prf}
By theorem \ref{th:comult}, $\delta$ is a map of left $\Ae$-modules. Since
$\delta$ has a counit, it is certainly injective.
\end{prf}

\section{Modules and Comodules}

We let $1_A \colon K \rightarrow A$ denote the map sending $1_K$ to
$1_A$. Since $X$ and $X \otimes K$ are canonically isomorphic, for any
map $f \colon X \rightarrow X$ we will abuse notation and write $f
\otimes 1_A \, \colon X \rightarrow X \otimes A$ instead of $f \otimes
1_A \, \colon X \otimes K \rightarrow X \otimes A$. When discussing
compositions of maps, the term ``switch'' will always refer to
reversing the order of noninteracting maps.

Suppose $M$ is a right $A$-module with structure map $m \colon M
\otimes A \rightarrow M$. Define the map $\comod_m \colon M
\rightarrow M \otimes A$ to be the composition:
\[
M    \myrightarrow{\cdot \otimes 1_A}
M \otimes A    \myrightarrow{\cdot \otimes \delta}
M \otimes A^2    \myrightarrow{m \otimes \cdot}
M \otimes A
\]

\begin{lemma}
   \label{lem:comodmod} The map $\comod_m$ endows $M$ with the
structure of a right $A$-comodule.
\end{lemma}
\begin{prf}
It is necessary to show that the following diagram commutes:
\begin{equation}
   \label{diag:comod} 
\begin{diagram}
   \node{M} \arrow{e,t}{\comod_m} \arrow{s,r}{\comod_m}
      \node{M \otimes A} \arrow{s,r}{\cdot \otimes \delta} \\
   \node{M \otimes A} \arrow{e,t}{\comod_m}
      \node{M \otimes A^2}
\end{diagram}
\end{equation}
Expanding each of the occurrences of $\comod_m$ in accordance with the
definition of that map yields the outer edge of this diagram:
\[
\begin{diagram}
   \node{M}  \arrow{e,t}{\cdot \otimes 1_A}
             \arrow{s,r}{\cdot \otimes 1_A}
      \node{M \otimes A} \arrow{e,t}{\cdot \otimes \delta}
                         \arrow{s,r}{\cdot \otimes \delta}
      \node{M \otimes A^2} \arrow{e,t}{m \otimes \cdot}
                      \arrow{s,r}{\cdot \cdot \otimes \delta}
      \node{M \otimes A} 
                      \arrow{s,r}{\cdot \otimes \delta}
   \\
   \node{M \otimes A}  \arrow{e,t}{\cdot \otimes \delta} 
                       \arrow{s,r}{\cdot \otimes \delta}
      \node{M \otimes A^2}
                      \arrow{e,t}{\cdot \otimes \delta \otimes \cdot}
      \node{M \otimes A^3}
                      \arrow{e,t}{m \otimes \cdot \cdot}
      \node{M \otimes A^2}
   \\
   \node{M \otimes A^2}
         \arrow{e,t}{\cdot \cdot \otimes 1_A \otimes \cdot}
         \arrow{se,b}{m \otimes \cdot}  
      \node{M \otimes A^3}
         \arrow{e,t}{\cdot \otimes A \otimes \delta \otimes \cdot}
         \arrow{n,r}{\cdot \otimes \mu \otimes \cdot}
      \node{M \otimes A^4}
         \arrow{n,r}{\cdot \otimes \mu \otimes A \otimes \cdot}
         \arrow{e,t}{m \otimes \cdot \cdot \cdot}
      \node{M \otimes A^3}
         \arrow{n,r}{m \otimes \cdot \cdot}
   \\
   \node[2]{M \otimes A}
            \arrow{e,b}{\cdot \otimes 1_A \otimes \cdot}
      \node{M \otimes A^2}
            \arrow{ne,b}{\cdot \otimes \delta \otimes \cdot} 
\end{diagram}
\]
From left to right and top down, the squares inside this large diagram
commute for the following reasons: Vacuity, coassociativity of
$\delta$, switch, property of the multiplicative identity, $\delta$
being a module map, $m$ being a module map. The hexagon on the bottom
is commutative because it only involves a switch.

It follows that the outer edge forms a commutative square, \ie diagram
(\ref{diag:comod}) is commutative.
\end{prf}

Suppose now that $M$ is a right $A$-comodule, with comodule structure
map $\comod \colon M \rightarrow M \otimes A$. Define the map
$m_{\comod} \colon M \otimes A \rightarrow M$ to be the composition:
\[
M \otimes A  \myrightarrow{\comod \otimes \cdot}
M \otimes A^2    \myrightarrow{\cdot \otimes \mu}
M \otimes A   \myrightarrow{\cdot \otimes \epsilon}
M
\]

\begin{lemma}
   \label{lem:modcomod}
The map $m_{\comod}$ endows $M$ with the structure of a right
$A$-module.
\end{lemma}
\begin{prf}
It is necessary to show that the following diagram commutes:
\begin{equation}
   \label{diag:mod}
\begin{diagram}
   \node{M \otimes A^2}
         \arrow{e,t}{m_{\comod} \otimes \cdot}
         \arrow{s,r}{\cdot \otimes \mu}
      \node{M \otimes A}
         \arrow{s,r}{m_{\comod}}
   \\
   \node{M \otimes A}
         \arrow{e,t}{m_{\comod}}
      \node{M}
\end{diagram}
\end{equation}

Expanding each occurrence of $m_{\comod}$ in accordance with the
definition of that map yields the outer edge of the following diagram:

\[
\begin{diagram}
   \node{M \otimes A^2}
         \arrow{e,t}{\comod \otimes \cdot \cdot}
         \arrow{se,t}{\comod \otimes \cdot \cdot}
         \arrow[3]{s,r}{\cdot \otimes \mu}
      \node{M \otimes A^3}
         \arrow{e,t}{\cdot \otimes \mu \otimes \cdot}
         \arrow{se,t}{\comod \otimes \cdot \cdot \cdot}
      \node{M \otimes A^2}
         \arrow{e,t}{\cdot \otimes \epsilon \otimes \cdot}
         \arrow{se,t}{\comod \otimes \cdot \cdot}
      \node{M \otimes A}
         \arrow{se,t}{\comod \otimes \cdot}   \\
   \node[2]{M \otimes A^3}
         \arrow{e,t}{\cdot \otimes \delta \otimes  A \otimes \cdot}
         \arrow{se,t}{\cdot \otimes \mu \otimes \cdot}
         \arrow[2]{s,r}{\cdot \cdot \otimes \mu}  
      \node{M \otimes A^4}
         \arrow{e,t}{\cdot \otimes A \otimes \mu \otimes \cdot}
      \node{M \otimes A^3}
         \arrow{e,t}{\cdot \cdot \otimes \epsilon \otimes \cdot}
      \node{M \otimes A^2} 
         \arrow{s,r}{\cdot \otimes \mu}  \\
   \node[3]{M \otimes A^2}
         \arrow{ne,t}{\cdot \otimes \delta \otimes \cdot}
         \arrow[2]{e,t}{\cdot \otimes \mu}
      \node[2]{M \otimes A}
         \arrow{s,r}{\cdot \otimes \epsilon}\\
   \node{M \otimes A}
         \arrow{e,t}{\comod \otimes \cdot}
      \node{M \otimes A^2}
         \arrow[2]{e,t}{\cdot \otimes \mu}
      \node[2]{M \otimes A}
          \arrow{ne,t}{\cdot \cdot}
          \arrow{e,t}{\cdot \otimes \epsilon}
      \node{M}
\end{diagram}
\]
The subdiagrams of this diagram are commutative for the following
reasons: In the top row of squares, the leftmost square expresses the
comodule property of $\comod$. The other two squares simply involve
switches, as does the large square on the far left. The square in the
center (between the second and third rows of maps) uses the module
property of $\delta$. The square to its right uses the counit property
of $\epsilon$. The large pentagon on the bottom expresses the
associativity of $\mu$. The triangle in the lower right hand corner is
vacuous.

It follows that the outer edge forms a commutative square, \ie diagram
(\ref{diag:mod}) is commutative.
\end{prf}

Lemmas \ref{lem:comodmod} and \ref{lem:modcomod} show that there are
canonical maps between the category of modules over $A$ and the
category of comodules over $A$. In fact, these provide an isomorphism.

\begin{thm}
   \label{th:modcomod} 
The category of right modules over a Frobenius algebra $A$ is
isomorphic to the category of right comodules over $A$.
\end{thm}
\begin{prf}
First we will show that the constructions $m \mapsto \comod_m$ and
$\comod \mapsto m_{\comod}$ are mutual inverses. Then we will show
that every module map is a comodule map for the corresponding comodule
structures, and vice-versa.

Suppose $m \colon M \otimes A \rightarrow M$ is a right module
structure map. Consider the following diagram:

\[
\begin{diagram}
   \node{M \otimes A}
         \arrow{e,t}{\cdot \otimes 1_A \otimes \cdot}
         \arrow{s,r}{\cdot \otimes \delta}
      \node{M \otimes A^2}
         \arrow{e,t}{\cdot \otimes \delta \otimes \cdot}
         \arrow{s,r}{\cdot \cdot \otimes \delta}
      \node{M \otimes A^3}
         \arrow{e,t}{m \otimes \cdot \cdot}
         \arrow{s,r}{\cdot \cdot \cdot \otimes \delta}
      \node{M \otimes A^2}
         \arrow{se,t}{\cdot \otimes \mu}
         \arrow{s,r}{\cdot \cdot \otimes \delta}
   \\
   \node{M \otimes A^2}
         \arrow{e,t}{\cdot \otimes 1_A \otimes \cdot \cdot}
         \arrow{s,r}{\cdot \cdot \cdot}
      \node{M \otimes A^3}
         \arrow{e,t}{\cdot \otimes \delta \otimes A \otimes \cdot}
         \arrow{sw,b}{\cdot \otimes \mu \otimes \cdot}
      \node{M \otimes A^4}
         \arrow{e,t}{m \otimes \cdot \cdot \cdot}
         \arrow{s,l}{\cdot \otimes A \otimes \mu \otimes \cdot}
      \node{M \otimes A^3}
         \arrow{s,l}{\cdot \otimes \mu \otimes \cdot}
      \node{M \otimes A}
         \arrow[2]{s,r}{\cdot \otimes \epsilon}
   \\
   \node{M \otimes A^2}
          \arrow[2]{e,b}{\cdot \otimes \delta \otimes \cdot}
      \node[2]{M \otimes A^3}
          \arrow{e,b}{m \otimes \cdot \cdot}
          \arrow{se,b}{\cdot \cdot \otimes \epsilon \otimes \epsilon}
      \node{M \otimes A^2}
          \arrow{ne,b}{\cdot \otimes \epsilon \otimes \cdot}
    \\
    \node[4]{M \otimes A}
          \arrow{e,t}{m}
      \node{M}
\end{diagram}
\]

\noindent
The composition of maps across the top and down the right is nothing
other than the definition of the map $m_{\comod_m} \colon M \otimes A
\rightarrow M$. Since the composition of maps down the left and across
the bottom is $m$ itself (by the counit property), the identity
$m_{\comod_m} \equiv m$ will follow if the diagram is
commutative. This is in fact the case, because the subdiagrams are
commutative for the following reasons: With the exception of those
that will now be mentioned explicitly, the subdiagrams are
commutative simply because they involve switches. The triangle on the
lower left uses the multiplicative unit property. The square to its
right expresses the module property of $\delta$. The square on the far
upper right is commutative because it is essentially the outer edge of
the following diagram:
\[
\begin{diagram}
   \addtolength{\dgARROWLENGTH}{-1.4em}
   \addtolength{\dgVERTPAD}{-1ex}
   \node[2]{A^2}
         \arrow{sw,t}{\cdot \otimes \delta}
         \arrow{s,r}{\mu}
         \arrow{se,t}{\mu}
   \\ 
   \node{A^3}
         \arrow{se,b}{\mu \otimes \cdot}
      \node{A}
         \arrow{s,r}{\delta}
      \node{A}
   \\
   \node[2]{A^2}
         \arrow{ne,b}{\epsilon \otimes \cdot}
\end{diagram}
\]
This latter diagram is commutative because the square on the left
expresses the module property of $\delta$, and the square on the right
express the counit property of $\epsilon$.

It follows that $m_{\comod_m} \equiv m$. Suppose, on the other hand,
that $\comod \colon M \rightarrow M \otimes A$ is a comodule
structure. We now show that $\comod_{m_{\comod}} \equiv
\comod$. Consider the following diagram:
\[
\begin{diagram}
   \node{M}
         \arrow{e,t}{\cdot \otimes 1_A}
         \arrow{se,b}{\comod} 
      \node{M \otimes A}
         \arrow{e,t}{\cdot \otimes \delta}
         \arrow{se,b}{\comod}
      \node{M \otimes A^2}
         \arrow{e,t}{\comod \otimes \cdot \cdot}
      \node{M \otimes A^3}
         \arrow{e,t}{\cdot \otimes \mu \otimes A}
      \node{M \otimes A^2}
         \arrow{s,r}{\cdot \otimes \epsilon \otimes \cdot}
   \\
   \node[2]{M \otimes A}
         \arrow{e,b}{\cdot \otimes \cdot \otimes 1_A}
      \node{M \otimes A^2}
         \arrow{ne,b}{\cdot \otimes A \otimes \delta}
         \arrow{e,b}{\cdot \otimes \mu}
      \node{M \otimes A}
          \arrow{ne,b}{\cdot \otimes \delta}
          \arrow{e,b}{\cdot \cdot}
      \node{M \otimes A}
\end{diagram}
\]
From left to right, the subdiagrams are commutative for the following
reasons: Switch, switch, the module property of $\delta$, the counit
property of $\epsilon$. Because the composition of maps across the top
and down the right of this diagram is simply the definition of
$\comod_{m_{\comod}}$, and the composition of maps down the left and
across the bottom is just $\comod$ (by the unit property of $1_A$), we
see that $\comod_{m_{\comod}} \equiv \comod$.

Suppose that $M$ and $N$ are right $A$-modules with module structure
maps $m$ and $n$ respectively. In order to verify that a map $f \colon
M \rightarrow N$ of right modules is also a map of right comodules
(for the corresponding comodule structures), consider the following
diagram:
\[
\begin{diagram}
   \node{M}
      \arrow{e,t}{\cdot \otimes 1_A}
      \arrow{s,r}{f}
   \node{M \otimes A}
      \arrow{e,t}{\cdot \otimes \delta}
      \arrow{s,r}{f \otimes \cdot}
   \node{M \otimes A^2}
      \arrow{e,t}{m \otimes \cdot}
      \arrow{s,r}{f \otimes \cdot \cdot} 
   \node{M \otimes A}
     \arrow{s,r}{f \otimes \cdot}
  \\
   \node{N}
      \arrow{e,t}{\cdot \otimes 1_B}
   \node{N \otimes A}
      \arrow{e,t}{\cdot \otimes \delta}
   \node{N \otimes A^2}
      \arrow{e,t}{n \otimes \cdot}
   \node{N \otimes A}
\end{diagram}
\]   
Two of the subdiagrams simply involve switches. The third is
commutative because $f$ is a map of modules. Thus, the outer edges
form a commutative diagram as well. But this diagram asserts that $f$
is a map of comodules, where the comodule structure maps are
$\comod_m$ and $\comod_n$.

If $f \colon M \rightarrow N$ is assumed to be a map of right
comodules, where the comodule structure maps are $\comod$ and
$\comod'$, then by reasoning analogous to that of the previous
paragraph, the following diagram shows that $f$ is a map of right
modules:
\[
\begin{diagram}
   \node{M \otimes A}
      \arrow{e,t}{\comod \otimes \cdot}
      \arrow{s,r}{f \otimes \cdot}
   \node{M \otimes A^2}
      \arrow{e,t}{\cdot \otimes \mu}
      \arrow{s,r}{f \otimes \cdot \cdot }
   \node{M \otimes A}
      \arrow{e,t}{\cdot \otimes \epsilon}
      \arrow{s,r}{f \otimes \cdot} 
   \node{M}
     \arrow{s,r}{f}
  \\
   \node{N \otimes A}
      \arrow{e,t}{\comod' \otimes \cdot}
   \node{N \otimes A^2}
      \arrow{e,t}{\cdot \otimes \mu}
   \node{N \otimes A}
      \arrow{e,t}{\cdot \otimes \epsilon}
   \node{N}
\end{diagram}
\]   
This completes the proof.
\end{prf}

With appropriate changes, all the results and proofs in this section
apply to left modules and left comodules as well.

\section{Cotensor Product}
   \label{sect:cotensor}

Suppose that $M$ is a right $A$-module with module structure map $m$,
and that $N$ is a left $A$-module with module structure map
$n$. By theorem \ref{th:modcomod}, $M$ is a right comodule with
structure map $\comod_m$ and $N$ is a left comodule with structure map
$\comod_n$. Let $\phi$ denote the map
\[
   \phi := \ \comod_m \! \otimes N \, - \, M \otimes \! \comod_n
        \, \colon \, M \otimes N  
        \longrightarrow
        M \otimes A \otimes N .
\]
The {\bf cotensor product} \cite{EilMoo66} $M \Box N$ of $M$ and $N$
is defined to be the kernel of $\phi$.

Viewing $A$ as both the right and left regular modules over itself
(\ie the module structure maps are both $\mu$), we can form $A \Box
A$. Note that $\comod_{\mu}$ is just the map $\delta$, by the module
property of $\delta$.

\begin{prop}
   \label{pr:ABoxA}
The cotensor product $A \Box A$ is exactly $\delta(A)$.
\end{prop}
\begin{prf}
By the definition of $\phi$, to show that $\delta(A) \subseteq A \Box
A$ it suffices to show that the two maps $(\comod_{\mu} \otimes A)
\circ \delta$ and $(A \otimes \comod_{\mu}) \circ \delta$ are the
same. But these two maps are just $(\delta \otimes A) \circ \delta$
and $(A \otimes \delta) \circ \delta$, respectively. These are the
same, by the coassociativity of $\delta$.

Now consider any element $x := \sum_i a_i \otimes b_i \, \in A \Box
A$. We have $(\delta \otimes A)x = (A \otimes \delta)x$, and thus 
\[
x \, = \, (\epsilon \otimes A^{\otimes 2}) \circ (\delta \otimes A)x 
\, = \, (\epsilon \otimes A^{\otimes 2}) \circ (A \otimes \delta)x
\, = \, \sum_i \epsilon(a_i)\delta(b_i) \, .
\]
It follows that $A \Box A \subseteq \delta(A)$.
\end{prf}

\begin{definition}
Let $D$ denote the left $\Ae$-submodule of $A \otimes A$ generated by
$T \circ \delta(1_A)$.  Note that $D$ and $\delta(A)$ (see corollary
\ref{cor:delta} above) are different objects.
\end{definition}

For any Frobenius algebra $A$, the map $\eta \colon A \otimes A
\rightarrow K$ defined by $\eta(a \otimes b) := \lambdal(1_A)(ab)$ is
a nondegenerate associative bilinear form \cite{CurRei}. If $\eta
\equiv \eta \circ T$, then $A$ is called a {\bf symmetric algebra}
[{\it ibid}].

\begin{prop}
   \label{pr:symmdelta} If $A$ is a symmetric algebra, then $D$ and
$\delta(A)$ are the same left $\Ae$-module.
\end{prop}
\begin{prf}
It suffices to show that if $A$ is a symmetric algebra, then
$\delta(1_A)$ is symmetric, \ie $T \circ \delta(1_A) =
\delta(1_A)$. Let $e_1, \ldots , e_n$ denote a basis for $A$, and let
$e_1^{\#}, \ldots , e_n^{\#}$ denote the dual basis of $A$ relative to
$\eta$, \ie the basis satisfying $\eta(e_i^{\#} \otimes e_j) =
\delta_{ij}$. The proof of proposition 5 in \cite{Abr96} (bearing in
mind the adjustments made in the proof of theorem \ref{th:comult} here
for noncommutativity) shows that $\delta(1_A) \ = \ \sum_j e_j \otimes
e_j^{\#}$.  Since, by assumption, we have $\eta(e_i^{\#} \otimes e_j)
= \eta(e_j \otimes e_i^{\#})$, a change of basis shows that
$\delta(1_A) \ = \ \sum_i e_i^{\#} \otimes e_i = T \circ \delta(1_A)$.
\end{prf}

\begin{cor}
   \label{cor:Ddelta}
If $A$ is commutative or semisimple or a group algebra then $D$ and
$\delta(A)$ are the same left $\Ae$-module.
\end{cor}
\begin{prf}
If $A$ is commutative then it is surely a symmetric algebra. Thus the
hypothesis of proposition \ref{pr:symmdelta} is automatically
satisfied.

By Wedderburn's first structure theorem, to prove the result in the
case when $A$ is semisimple it suffices to assume that $A$ is a matrix
ring. In that case, $A$ has a Frobenius algebra structure given by the
map $\lambdal(1_A)(a) := \mathrm{Tr}(a)$. It is an easy exercise to
show that this provides $A$ with the structure of a symmetric algebra.

In the case of a group algebra $A$ over group $G$, the Frobenius
algebra structure is given by the map which returns the coefficient of
the identity element. The coproduct then sends $1_A$ to $\sum_{g \in
G}g \otimes g^{-1}$, which is clearly symmetric.
\end{prf}
 
When $A$ is not a symmetric algebra, proposition \ref{pr:symmdelta}
does not necessarily apply.

\begin{example}
   \label{ex:1}
Let $A$ denote the exterior algebra on two generators, $x$ and
$y$. Then
\[
\delta(1_A) \, = \, 1_A \otimes xy \, + \, xy \otimes 1_A 
          \, - \, x \otimes y \, + \, y \otimes x ,
\]
and $\delta(A)$ has the basis 
\[
   \left\{ \delta(1_A),\ \ 
      x \otimes xy \, + \, xy \otimes x,\ \, \ 
      y \otimes xy \, + \, xy \otimes y,\ \, \ 
      xy \otimes xy \right\} \, ,
\]
whereas $D$ has the basis
\[
   \left\{ T \circ \delta(1_A),\ \ 
      x \otimes xy \, - \, xy \otimes x,\ \, \
      y \otimes xy \, - \, xy \otimes y,\ \, \
      xy \otimes xy \right\} \, .
\]
\end{example}
     
\vspace{3pt}

Given a right module $M$ and a left module $N$ as above, endow $N
\otimes M$ with the obvious left $\Ae$-module structure. Let $\Hom{D}$
denote the vector space of left $\Ae$-module maps.

\begin{thm}
   \label{th:cotensorhom} 
There is a vector space isomorphism
\[
M \Box N \ \cong \ \Hom{D} \, .
\]
\end{thm}
\begin{prf}
Note first that an element $f \in \Hom{D}$ is determined by its value
on $T \circ \delta(1_A)$, the generator of $D$. 

The following diagram is commutative, since $f$ is a map of modules:
\[
\begin{diagram}
   \addtolength{\dgARROWLENGTH}{-1.4em}
    \node{K} 
         \arrow[2]{e,b}{\delta(1_A) \, \otimes \, T \circ \delta(1_A)}
      \node[2]{A^4} 
         \arrow[3]{e,b}{(\cdot \cdot \otimes \mu) \circ T_{1432} \ 
                   - \ (\cdot \otimes \mu \otimes \cdot)}
         \arrow{s,r}{\cdot \cdot \otimes f}
      \node[3]{A^3}
         \arrow{s,r}{\cdot \otimes f}
   \\
   \node[3]{A^2 \otimes N \otimes M}
         \arrow[3]{e,t}{(\cdot \cdot \otimes m) \circ T_{1432} \
                   - \ (\cdot \otimes n \otimes \cdot)}
      \node[3]{A \otimes N \otimes M}
\end{diagram}
\]
By the comodule property of $\delta$, the composition of maps across
the top of the diagram is $0$. Since the composition of maps from the
upper left, down and across the bottom is $T_{132} \circ \phi \circ f
\left[T \circ \delta(1_A) \right] $, it follows that $f \left[ T \circ
\delta(1_A) \right] \in M \Box N$. Thus, there is a well defined
injective map $\sigma \colon \Hom{D} \rightarrow M \Box N$ sending $f
\mapsto f \left[ T \circ \delta(1_A) \right]$. Since each element $e
\in N \otimes M$ defines a unique $\Ae$-module map $\tau(e) \colon T
\circ \delta(1_A) \mapsto e$, restriction of $\tau$ to $M \Box N$
provides an inverse to $\sigma$.
\end{prf}

Allowing for abuse of notation, define the cotensor product functor
$\Box_A$ by $\Box_A \colon M \otimes N \mapsto M \Box N$, and let
$\mathrm{Cotor}^i_A(M,N)$ denote its right derived functors. Let
$H^i(A, -)$ denote the Hochschild cohomology functors.

\begin{cor}
   \label{cor:cotorext}
Over a Frobenius algebra $A$, the Cotor functor is given by
\[
\mathrm{Cotor}^*_A(M,N) \cong \mathrm{Ext}^*_{A^e}(D,N \otimes M) \, .
\]
\end{cor}
\begin{prf}
In light of theorem \ref{th:cotensorhom}, this is purely a matter of
definitions.
\end{prf}

\begin{cor}
   \label{cor:cotorHoch} 
If $A$ is a symmetric algebra, then cotensor product and its derived
functors are Hochschild cohomology, \ie
\[
\mathrm{Cotor}^*_A(M,N) \cong H^*(A,N \otimes M)\, . 
\]
\end{cor}
\begin{prf}
By corollary \ref{cor:Ddelta} we have $\delta(1_A) = T \circ
\delta(1_A)$ and thus $D = \delta(A)$. Since, by corollary
\ref{cor:delta}, $\delta$ is an injective map of left $\Ae$-modules
(determined by its value on $\delta(1_A)$), $D$ and $A$ are isomorphic
as $\Ae$-modules. It follows from theorem \ref{th:cotensorhom} that $M
\Box N \cong \Hom{A}$. But this is exactly $H^0(A,N \otimes M)$
\cite[pg. 301]{Wei}. Since $H^*(A,-) \cong \mathrm{Ext}^*_{A^e}(A,-)$
[{\it ibid.} pg. 303], this corollary follows from \ref{cor:cotorext}.
\end{prf}


\end{document}